\documentclass[11pt]{article} 

\usepackage[utf8]{inputenc} 

\usepackage{geometry} 
\usepackage{graphicx} 

\usepackage{booktabs} 
\usepackage{array} 
\usepackage{paralist} 
\usepackage{verbatim} 
\usepackage{subfig} 
\usepackage{amsmath}
\usepackage{amssymb}
\usepackage{nicematrix}

\usepackage{fancyhdr} 
\usepackage{sectsty}
\allsectionsfont{\sffamily\mdseries\upshape} 

\usepackage[nottoc,notlof,notlot]{tocbibind} 
\usepackage[titles,subfigure]{tocloft} 

\title{Some NP Complete Problems Based on Algebra and Algebraic Geometry}
\author{Paul Hriljac\\ Department of Mathematics, ERAU}
\date{} 

\begin{document}
\maketitle

\noindent 
\textbf{Notation}: Let  $p$ be a prime and let $\mathbb{F} $ be a finite field with $q=p^n \geq 4$ elements and algebraic closure $\mathbb{\bar{F}} $. For sets $X,Y$ let $Hom(X,Y)$ denote the functions from $X$ to $Y$. For numbers $m,n$ let $Hom(m,n)$ denote $Hom(\{1,\dots m\},\{1,\dots, n\})$. Let $\Gamma$ be a connected graph with vertex set $V$ and edge set $E$.
\newline
\newline
\textbf{Definition} a $\mathbb{F}$\textit{-coloring} of $\Gamma$ is an assignment  \(V \rightarrow \mathbb{F^{*}}\) so that no adjacent vertices share the same value.
\newline
\newline
\noindent From now on, suppose that each edge of $\Gamma$ has be given an orientation, so one can speak of the initial vertex and terminal vertex of each edge.

\noindent Assign to $\Gamma$ the array $A_{\Gamma} \in \mathbb{F}^{|V||E|(q-1)}$ given by 
\begin{equation*}
A_{\Gamma}(v,e,c)=
\begin{cases} 
 &1  \text{ if }  v=e(0)\\
& c \text{ if } v=e(1)\\
& 0 \text{ else}.
\end{cases}
\end{equation*}
Here $v$  is vertex, $e $ is an edge and $c$ is an element of $\mathbb{F}-\{0\}$.
Let $\mathbf{D}_{V,E}$ be the determintal variety corresponding to $|V|$ by $|E|$ matrices with rank $<min(|E|,|V|)$ (see [3]). We suppress the indices when the context is clear. Then $\mathbf{D}$ is a subscheme of $\mathbb{A}^{|V||E|}$. For $\alpha \in Hom(E,\mathbb{F}-\{0,-1\})$ let $\pi_{\alpha}$ be the projection $\mathbb{A}^{|V||E|(q-1)} \rightarrow \mathbb{A}^{|V||E|}$ given by $A \mapsto [A(v,e,\alpha(e))]_{v,e}$. Let $\mathbf{D}_{\alpha}=\pi_{\alpha}^{*}\mathbf{D}$ and let $\mathbf{C}=\mathbf{C}_{|V||E|(q-1)}=\cup_{\alpha}\mathbf{D}_{\alpha}$. Call $\mathbf{C}$ the \textit{coloring variety} associated to graphs with $|E|$ edges, $|V|$ vertices and $q-1$ colors. Since $\mathbf{C}$ can be defined by cofactor relations, it can be realized as a closed subscheme of $\mathbb{A}^{|V||E|(q-1)}$.
\newline
\newline
\textbf{Theorem}: \textit{Assuming} $|V| \leq |E|$ $\Gamma$ \textit{is}  $\mathbb{F}$-\textit{colorable if and only if} $A_{\Gamma} \in \mathbf{C}$.
\newline
\newline \textbf{Proof}: Given an $\mathbb{F}$\textit{-coloring} $\{c_v\}$ of $\Gamma$, Fermat's Theorem implies $c_{v}^{q-1}=1$ for all vertices $v$. Therefore, if $v$ and $w$ are adjacent 
\begin{equation*}\frac{c_{v}^{q-1}-c_{w}^{q-1}}{c_v-c_w}=0. \end{equation*}
Expanding this out gives \begin{equation*}\prod_{\alpha \in \mathbb{F}-\{0,-1\}} c_{e(0)}+\alpha c_{e(1)}=0. \end{equation*}
Therefore $c_{e(0)}+\alpha c_{e(1)}=0$ for some $\alpha =\alpha(e)\in \mathbb{F}-\{0,-1\}$, which implies $\mathbf{c}M_{\Gamma}=0$ where
\begin{equation*}
M_{\Gamma}=\pi_{\alpha}(A_{\Gamma})=\begin{cases} 
 &1  \text{ if }  v=e(0)\\
& \alpha(e) \text{ if } v=e(1)\\
& 0 \text{ else}
\end{cases}, \qquad
\mathbf{c}=[c_v]_{v \in V}.
\end{equation*}
This shows that the rank of $\pi_{\alpha}(A_{\Gamma})$ is less than $|V|$ and hence $A_{\Gamma} \in \mathbf{C}$.
\newline
On the other hand if $\pi_{\alpha}(A_{\Gamma})$ has rank less than $|V|$ for some $\alpha \in Hom(E,\mathbb{F}-\{0,-1\})$ than a nonzero vector in the kernel space $\pi_{\alpha}(A_{\Gamma})$ is easily seen to provide a procedure yielding a valid coloring for the graph.
\newline
\newline
From this proof one sees the
\newline
\newline
\textbf{Theorem}: \textit{If $|E|<|V|$ then $\Gamma$ \textit{is}  $\mathbb{F}$-\textit{colorable}}.
\newline
\newline
\textbf{Proof}: Since $|E|<|V|$ the kernel space of the linear map corresponding to $\mathbf{x} \mapsto \mathbf{x}\pi_{\alpha}(A_{\Gamma})$ is nonzero for any $\alpha$. Using a nonzero element of the kernel space one can create a coloring.
\newline
\newline
\textbf{Example 1}: Try to 3-color the following graph:
\begin{figure}[htp]
    \centering
    \includegraphics[width=12cm]{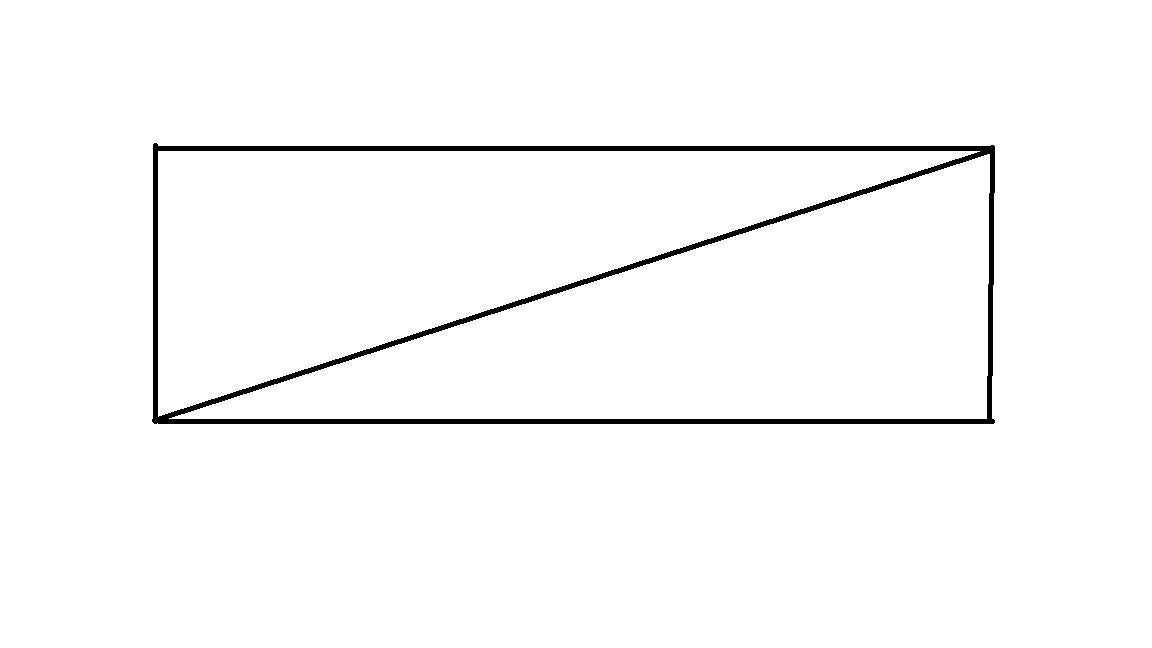}
\end{figure}
\newline The theorem implies that the graph has a 3-coloring iff there is a map $\alpha:\{1,2,3,4,5\} \rightarrow \mathbb{F}_{4}-\{0,-1\}$ such that the rank of the matrix
\begin{equation*}
\pi_{\alpha}(A_{\Gamma})= 
\begin{bmatrix}
 \alpha_{1} &1 &0 &0 &0\\
0 &\alpha_{2} &1 &0 &\alpha_{5}\\
1&0&0&\alpha_{4}&1\\
0&0&\alpha_{3}&1&0 
\end{bmatrix}.
\end{equation*}
is less than 4. Since the determinant of the matrix with the last column deleted is 1, the rank is never less than 4 so the graph cannot be 3-colored.
\newline
\newline
\textbf{Definition}: Given a family $\mathfrak{X}=(X,Y)$ of varieties $X$ and subvarieities $Y \subset X$ over a field $K$ \textit{the membership problem for} $\mathfrak{X}$ to be the following: Given $P \in X(K)$, is $P \in Y(K)$?
\newline
\newline
\textbf{Definition}: Let $\mathfrak{C}=(\mathbf{C}_{m,n,q-1},\mathbb{A}^{mn(q-1)})$ using coloring varities as defined above.
\newline
\newline
\textbf{Theorem}: \textit{The membership problem for} $\mathfrak{C}$ \textit{over a field} $\mathbb{F}$ \textit{with} $q \geq 4$ \textit{elements is NP Complete}.
\newline
\newline
\textbf{Proof}: To see that this problem is in NP, one just needs to check the rank of $\pi_{\alpha}(A_{\Gamma})$ for a given $\alpha$. To see that this problem is NP hard note that the NP complete problem of $\mathbb{F}$-coloring graphs is subsumed  by this problem.
\newline
\newline
We now introduce another NP complete membership problem with the following example: 
\newline
\newline
\textbf{Example 2}: Consider the graph in Example 1. Parametrize a line in $\mathbb{A}^{4,5,2} $ by 
\begin{equation*}
\gamma(t)=
\begin{bmatrix}
 \alpha_{1} &t &0 &0 &1-t\\
0 &\alpha_{2}t &1 &0 &\alpha_{5}\\
t&0&0&\alpha_{4}&t\\
0&0&\alpha_{3}&t&0 
\end{bmatrix}.
\end{equation*}
then $\gamma(1)=A_{\Gamma}$. Each choice of $\alpha$ yields a 4 by 5 array with coefficients in $\mathbb{F}[t]$. Taking the five major cofactors of this array gives five polynomials $f_{1,\alpha},\dots,f_{5,\alpha}$. Taking the greatest common divisor of these polynomials gives a polynomial $f_{\alpha}$ of degree $\leq 4$. It is easy to see 
\begin{equation*}
f_{\alpha}(1)=0 \text{ if and only if } \alpha \text{ gives a valid coloring of the graph.}
\end{equation*}
\newline
Next suppose $\mathcal{F}=\{f(x)_{\alpha}\}_{\alpha \in Hom(n,l)}$ is a set of polynomials in $\mathbb{F}[x]$ with degree $\leq n$ whose coefficients can be found from $\alpha$ in time polynomial in $l,m,n$. Let $Z(\mathcal{F})$ denote the roots of $\prod_{\alpha}f_{\alpha}$ and  $\mathfrak{Z}=(Z(\mathcal{F}),\mathbb{A}^1)$.
\newline
\newline
\textbf{Theorem}: \textit{The membership problem for} $\mathfrak{Z}$ \textit{over a field} $\mathbb{F}$ \textit{with} $q \geq 4$ \textit{elements is NP Complete}.
\newline
\newline
\textbf{Proof}: Using elementary linear algebra it is easy to check that the variety $D_{n,m}$ is cut out set-theoretically by the equations $\delta_{0}=\dots=\delta_{m-n}=0$ where $\delta_{i}$ is the determinant of the $n \times n$ matrix consisting of the first $n-1$ columns of a matrix as well as the column in position $n+1$. To decide if some $A \in \mathbb{A}^{m,n,l}$ is in $\mathbf{C}_{m,n,l}$ let $L$ be a line in $ \mathbb{A}^{m,n,l}$ parametrized by some linear function $\theta(t)$ with $\theta(t_0)=A$. Let $f_{\alpha}(t)=gcd[\delta_{0}(\pi_{\alpha}(\theta(t))),\dots,\delta_{m-n}(\pi_{\alpha}(\theta(t)))]$.  Then $A \in \mathbf{D}_{\alpha}$ iff $f_{\alpha}(t_0)=0$ and so $A \in \mathbf{C}$ iff $t_0$ is a root of $\prod_{\alpha}f_{\alpha}$. This shows that the $\mathfrak{Z}$ membership problem is NP-hard. It is easy to see that the problem is in NP: Given $\alpha$ and $t_0$, just check that $f_{\alpha}(t_0)=0$.
\newline
\newline
\textbf{Definition}: A \textit{computation tree} over a finite field $\mathbb{F}$ is a rooted tree binary with a collection of field elements associated with each node (see [1]). Other than roots and leaves, each node is either a computation or a branch. Each computation node performs an algebraic calculation over that field that depends only on the field elements associated with earlier nodes. The branch nodes have a single associated field element and the two edges emanating from it are  labeled $"="$ and $"\neq"$.  A \textit{decision tree} $T$ for $X \subset \mathbb{F}$ is a computation tree that decides if some $x \in \mathbb{F}$ belongs in $X$. The leaves are labeled $"\in"$ and $"\notin"$. Since the tree is operating only on elements of $\mathbb{F}$ we may assume that the root node has a unique field element.  Let $|T|$ denote the length of the longest path through $T$. 
\newline
\newline
\textbf{Definition}: A \textit{computation sequence} for a polynomial $f(x) \in \mathbb{F}[x]$ is a sequence of polynomials  $g_i(x) \in \mathbb{F}[x], i = 0,\dots,r$ such that:
\begin{align*}
&g_0=x;\\
&g_{i}=(\sum_{j=0}^{i-1}a_{i,j}g_j)*(\sum_{j=0}^{i-1}b_{i,j}g_j), \quad a_{i,j}, b_{i,j} \in \mathbb{F};\\
&f(x)=f(0)+\sum_{i=0}^{r}c_i*g_i, \quad c_i \in \mathbb{F};
\end{align*}
Let $\L(f) $ denote the length of the shortest computation sequence for $f$. This is the \textit{multiplicative complexity} of a polynomial $f(x) \in \mathbb{F}[x]$.
\newline
\newline
Suppose $\mathbb{F}$ is a finite field with $q$ elements and $X \subset \mathbb{F}$ a subset with  characteristic function $\chi$. Then $\chi$ can be realized as a polynomial $\chi(x) \in  \mathbb{F}[x]$. Suppose $T$ is a decision tree for $X \subset \mathbb{F}$ . Then
\newline
\newline
\textbf{Theorem}: $L(\chi) \leq 24*log(q)*|T|$.
\newline
\newline
To prove the theorem we must introduce some new ideas and discuss their properties. 
\newline
\newline
\textbf{Definition}: Given a decision tree for $X \subset \mathbb{F}$ let $\mathbf{p}_x$ denote the path $T$ takes when input $x \in \mathbb{F}$ is provided. Let $|\mathbf{p}_x|$ denote the total number of internal nodes in $\mathbf{p}_x$.  Let $\mathbf{p}_T$ be the \textit{generic path} of $T$, corresponding to the decision $"\neq"$ at all decision nodes of $T$.
\newline
\newline
\textbf{Definition}: Say that a decision tree $T$ is \textit{pruned} if :
\newline
1. For every decision node $\nu$ in $T$ there are inputs $x,y \in \mathbb{F}$ such that $\mathbf{p}_x,\mathbf{p}_y$ both pass through $\nu$ but take different forks there.
\newline
2. Every computation node has at least one decision node as a descendant.
\newline
\newline
\textbf{Lemma 1}:\textit{ Given any decision tree $T$ there is a pruned decision tree $T'$ deciding the same membership with $|T'| \leq |T|$}.
\newline
\newline
\textbf{Proof}: Simply remove the nodes of $T$ that are not needed.
\newline
\newline
\textbf{Definition}: A decision tree $T$ is \textit{polynomial} if no computation node performs division.
\newline
\newline
\textbf{Lemma 2}: \textit{Given a decision tree $T$ there is a pruned polynomial decision tree $T'$ deciding the same membership with $|T'| \leq 4|T|$}
\newline
\newline
\textbf{Proof}: If $T$ is not pruned, prune it, assume henceforth that it is pruned. If it has no nodes performing division, do nothing. Otherwise proceed as follows:
\newline 
\begin{enumerate}
\item Replace the input node $\nu$ with assignment $u_{\nu} \leftarrow x$ with an input node $\nu_{1}$ with assignment $u_{\nu_1} \leftarrow x$ and and a 0-arity node $\nu_{2}$ with $u_{\nu_{2} }\leftarrow 1$.
\item Replace every computation node $\nu$ with two nodes or four nodes as follows:
\begin{enumerate}
\item If $\nu$ is a 0-arity node $u_{\nu}\leftarrow c$ replace $\nu$ with nodes $\nu_1,\nu_2$ and let $u_{\nu_1}\leftarrow c, u_{\nu_2} \leftarrow 1$.

\item If $\nu$ is a 1-arity node $u_{\nu}\leftarrow c*u_{\mu}$ replace $\nu$ with nodes $\nu_1,\nu_2$ and let $u_{\nu_1}\leftarrow c*u_{\mu_1},u_{\nu_2}\leftarrow u_{\mu_2}$.
\item If $\nu$ is a 2-arity node $u_{\nu} \leftarrow u_{\mu}/u_{\lambda}$ replace $\nu$ with nodes $\nu_1,\nu_2$ and let $u_{\nu_1} \leftarrow u_{\mu_1}*u_{\lambda_2},u_{\nu_2} \leftarrow u_{\mu_2}*u_{\lambda_1} $.

\item If $\nu$ is a 2-arity node $u_{\nu}\leftarrow u_{\mu}+u_{\lambda}$ replace $\nu$ with nodes $\nu_1,\dots \nu_4$ and let $u_{\nu_3} \leftarrow u_{\mu_1}*u_{\lambda_2},u_{\nu_4} \leftarrow u_{\mu_2}*u_{\lambda_1}, u_{\nu_1} \leftarrow u_{\nu_3}+u_{\nu_4},u_{\nu_2} \leftarrow u_{\nu_3}*u_{\nu_4}$.

\item If $\nu$ is a 2-arity node $u_{\nu}\leftarrow u_{\mu}-u_{\lambda}$ replace $\nu$ with nodes $\nu_1,\dots \nu_4$ and let $u_{\nu_3} \leftarrow u_{\mu_1}*u_{\lambda_2},u_{\nu_4}  \leftarrow u_{\mu_2}*u_{\lambda_1}, u_{\nu_1} \leftarrow u_{\nu_3}-u_{\nu_4},u_{\nu_2} \leftarrow u_{\nu_3}*u_{\nu_4}$.

\item If $\nu$ is a 2-arity node $u_{\nu}\leftarrow u_{\mu}*u_{\lambda}$ replace $\nu$ with nodes $\nu_1,\nu_2$ and let $u_{\nu_1}\leftarrow u_{\mu_1}*u_{\lambda_1},u_{\nu_2}  \leftarrow u_{\mu_2}*u_{\lambda_2} $.

\end{enumerate}
\item For a decision node $\nu$ making the query $"u_{\mu}=u_{\lambda}?"$ add two computation nodes $\nu_1,\nu_2$ with $u_{\nu_1}=u_{\mu_1}*u_{\lambda_2},u_{\nu_2}=u_{\mu_2}*u_{\lambda_1}$ and then replace the decision node with a new node making the query $"u_{\nu_1}=u_{\nu_2}?"$.
\item Connect all leaf nodes to the modifed decision nodes corresponding to the previous logic.
\end{enumerate}
This procedure reflects the fact that $u_{\nu}=\frac{u_{\nu_1}}{u_{\nu_2}}$ and all operations using $u_{\nu}$ can be written in terms of $u_{\nu_1}$ and $u_{\nu_2}$. Hence the value $T$ calculates at a computation node with an input $x \in \mathbb{F}$ corresponds exactly to the values $T'$ calculates at the corresponding nodes with the same input. Any decision node of $T'$ is making a comparison on the same values as the corresponding node of $T$. It's clear that if $T$ is pruned, so is $T'$. QED.
\newline
\newline
We introduce a new operation for modifying decision trees in a way that trades computation nodes for decsion nodes. Let $T$ be a pruned polynomial decision tree for deciding membership of $ X\subset\mathbb{F}$. Suppose $\lambda $ is a leaf node of maximal distance from $\mathbf{p}_T$, as measured in the number of decisions nodes between $\lambda$ and $\mathbf{p}_T$ and suppose $\lambda \notin \mathbf{p}_T$. Let $\mu$ the parent decision node of $\lambda$ and suppose $\mu\notin\mathbf{p}_T$. Let $\nu$ be the parent decision node of $\mu$. Since $T$ is pruned there are two cases:
\begin{enumerate}
\item One child of $\nu$ is $\mu$. One child of $\mu$ is $\lambda$. The other child of $\mu$ is another leaf node $\lambda 1$ with a label that is the opposite of the label of $\lambda$. The other child of $\nu$ is a leaf node $\lambda 2$. If necessary, switch the names of leaf nodes $ \lambda,\lambda 1$  so that $\lambda$ is the leaf node labeled $"\in"$.
\item One child of $\nu$ is $\mu$. The other child of $\nu$ is another decision node $\mu 1$. One child of $\mu$ is $\lambda$. The other child of $\mu$ is another leaf node $\lambda 1$ with a label which is the opposite of $\lambda$. The decision node $\mu 1$ has children, which are both leaf nodes and have opposite labels.
\end{enumerate}

\noindent We describe a procedure RETRACT that replaces the entire segment of the tree which emanates from $\nu$ with a segment which has only one decision node $\xi$ whose children are two leaf nodes. The new segment will have two leaves as children and will be labeled so that any input $x \in \mathbb{F}$  whose path enters that segment will now enter the new segment and terminate at a leaf with the same label as it had previously.
\newline
\newline
RETRACT Case I.  
The following table assumes that the leaf node $\lambda$ is always labeled $"\in"$ and its sibling is always labeled $"\notin"$. The label for leaf node varies as indicated. The middle column describes the path emanating from decision node $\nu$ and terminating at leaf node $\lambda$. The third column describes the decision function that the new node $\xi$ is using.
\newline
\begin{center}
\begin{tabular}{ c c c }
 $\lambda 2$ label & decisions for $\nu \rightarrow \lambda$ & decision for $\xi$ \\ 
 $\in$ & $u_\nu=0,u_\mu=0$ & $(1-u_\nu^{q-1})u_\mu$=0? \\  
 $\in$ & $ u_\nu=0,u_\mu\neq 0$ & $ (1-u_\nu^{q-1})(1-u_\mu^{q-1})=0?$   \\
 $\in$ & $u_\nu\neq0,u_\mu=0$ & $u_\mu u_\nu=0? $\\
$\in$ & $u_\nu\neq 0,u_\mu\neq 0 $& $(1-u_\mu^{q-1})u_\nu=0?$ \\
$\notin$ & $u_\nu=0,u_\mu=0$ & $1-(1-u_\nu^{q-1})(1-u_\mu^{q-1})=0?$\\
$\notin$ & $u_\nu=0,u_\mu\neq 0$ & $1-(1-u_\nu^{q-1})u_\mu^{q-1}=0? $\\
$\notin$ & $u_\nu\neq 0,u_\mu=0$ & $1-(1-u_\mu^{q-1})u_\nu^{q-1}=0?$\\
$\notin$ & $u_\nu\neq0,u_\mu\neq 0$ & $1-u_\nu^{q-1}u_\mu^{q-1}=0?$
\end{tabular}
\end{center}

\noindent RETRACT Case II
The following table describes the situation where the decision node has two children, both of which are decision nodes and each of these nodes in turn have two children, both of which are leaf nodes. The left column describes the labels of the four leaves numbered so that they correspond to the decisions first node is $0$, second node is $0$; first node is $0$, second node is $\neq 0$; first node is $\neq 0$, second node is $0$; first node is $\neq 0$, second node is $\neq 0$. The second column describes the decision function that the new node $\xi$ is using.

\begin{center}
\begin{tabular}{ c c }
labels on $\lambda,\lambda 1,\lambda 2,\lambda 3$ &decision for $\xi$ \\
$\notin \in \in \notin$&$(1-u_\nu^{q-1})u_\mu+u_\mu u_\nu=0?$\\
$\notin \in \in \notin$&$(1-u_\mu^{q-1})(1-u_\nu^{q-1})+u_\mu u_\nu=0?$\\
$\in \notin \notin \in$&$(1-u_\nu^{q-1})u_\mu+(1-u_\mu^{q-1})u_\nu=0?$\\
$\notin \in \notin \in$&$(1-u_\mu^{q-1})(1-u_\nu^{q-1})+(1-u_\mu^{q-1})u_\nu=0?$
\end{tabular}
\end{center}
In both cases the process of replacing the segment of $T$ that emanates from the node $\nu$ and replacing it by a node $\xi$ which has two children, both leaf nodes with opposite labels, is said to apply RETRACT on $T$ at $\nu$. 
\newline
\newline
For a decision node $\nu$, let $ |\nu| $ denote the number of decision nodes that separate $\nu$ from the generic path $\mathbf{p}_T$. Let $\langle T \rangle = max_{\nu}|\nu|$ where $\nu$ ranges over all decision nodes of $T$. 
\newline
\newline
\textbf{Lemma 3}: \textit{If $T$ is a  decision tree for $X \subset \mathbb{F}$  that has $\langle T \rangle >1$  then applying RETRACT on $T$ at all decision nodes $\nu$ with $|\nu|=\langle T \rangle $ results in a new tree $T'$ with $\langle T' \rangle=\langle T \rangle-1$ and $|T'| \leq |T|+6log(q)$.}
\newline
\newline
\textbf{Proof}: This is immediate from the construction.
\newline
\newline
\textbf{Definition}: Suppose  $T$ is a  decision tree for $X \subset \mathbb{F}$  with  $\langle T \rangle \leq 1$, so every leaf is a child of a decision node on $\mathbf{p}_T$ or a child of a decision node that is seperated from $\mathbf{p}_T$ only by computation nodes. Let $D_0$ be the set of decision nodes on $\mathbf{p}_T$. Let $D_1$ denote the decision nodes $\nu$ not on $\mathbf{p}_T$ such that when $x \in \mathbb{F}$ is entered $u_{\nu}=0$ implies $x \in X$. Let $D_2$ denote the decision nodes $\nu$ not on $\mathbf{p}_T$ such that when $x \in \mathbb{F}$ is entered $u_{\nu}\neq 0$ implies $x \in X$. For $\nu \in D_1 \cup D_2$ let $\mu(\nu) \in \mathbf{p}_T$ denote the decision node that is the parent node (disregarding computation nodes).
\newline
\newline
\textbf{Lemma 4}: \textit{Suppose  $T$ is a  decision tree for $X \subset \mathbb{F}$  with  $\langle T \rangle \leq 1$. Then}
\begin{equation*}
\chi = \prod_{\nu \in D_0}u_{\nu} \times \prod_{\nu \in D_1}(1-(1-u_{\nu}^{q-1})(1-u_{\mu(\nu)}^{q-1}))\times \prod_{\nu \in D_2}(1-u_{\nu}(1-u_{\mu(\nu)}^{q-1})
\end{equation*}
\textbf{Proof of Theorem}: First turn $T$ into a polynomial pruned tree. Then apply RETRACT until one obtains a tree $T'$ with  $\langle T' \rangle \leq 1$. Now apply Lemma 4 to obtain the characteristic function of $X \subset \mathbb{F}$. Due to Lemmas 1 to 4, this gives the required upper bound on the complexity of $\chi$. 
\newline
\newline
We constrast this last theorem with the following result of Patterson and Stackmeyer (see [2]).
\newline
\newline
\textbf{Theorem}: \textit{For each degree $d$ there exists a hypersurface $H_d$ in the space of  polynomials $f(x) \in \bar{\mathbb{F}}[x] $  of degree $d$ such that $f \notin H_d$ implies $L(f) \geq \sqrt
{d}-1$}
\newline
\newline
\textbf{Defintion}: We will call a polynomial $f(x) \in \bar{\mathbb{F}}[x] $ \textit{generic} if  $f \notin H_d$.
\newline
\newline
\textbf{Definition}: Given graphs $\Gamma_1, \Gamma_2$ of the same size with edge sets of cardinality $n$ and a finite field $\mathbb{F}$ with $q=p^n$ elements let $\theta(x)=A_{\Gamma_1}+x A_{\Gamma_2},$
$$
\Phi_{\Gamma_1,\Gamma_2}(x) =\prod_{\alpha}f_{\alpha}(x)
$$
where $\alpha: (1,\dots,n) \rightarrow \mathbb{F}-\{0,-1\}$ ranges over all possible maps and $f_{\alpha}(x)=
 gcd[\delta_{0}(\pi_{\alpha}(\theta(x))),$  $\dots,\delta_{n-v}(\pi_{\alpha}(\theta(x)))]$  and $\delta_{i}$ is the determinant of the $v \times v$ matrix consisting of the first $v-1$ columns of a matrix as well as the column in position $v+1$ 
Let $\chi_{\Gamma_1,\Gamma_2} ={\Phi_{\Gamma_1,\Gamma_2}}^{q-1}$.
\newline
\newline
\textbf{Conjecture}: \textit{For all $n$ there exisits graphs $\Gamma_1, \Gamma_2$ of the same size with edge sets of cardinality $n$ such that $\chi_{\Gamma_1,\Gamma_2}$ has degree exponential in $n$ and is generic.}
\newline
\newline
We observe that the conjecture immediately implies $NP \neq P$ (see [4]).
\newline
\newline
\textbf{References}:
\newline
1. Bürgisser, Peter; Clausen, Michael; Shokrollahi, M. Amin (1997). Algebraic complexity theory. Grundlehren der Mathematischen Wissenschaften. Vol. 315. With the collaboration of Thomas Lickteig. Berlin: Springer-Verlag. ISBN 978-3-540-60582-9. Zbl 1087.68568.
\newline
2. Paterson, Michael S.; Stockmeyer, Larry J. (1973). "On the Number of Nonscalar Multiplications Necessary to Evaluate Polynomials". SIAM Journal on Computing. 2 (1): 60–66. doi:10.1137/0202007.
\newline
3. Bruns, Winfried; Vetter, Udo (1988). Determinantal rings. Lecture Notes in Mathematics. Vol. 1327. Springer-Verlag. doi:10.1007/BFb0080378. ISBN 978-3-540-39274-3.
\newline
4.  Cook, Stephen (1971). "The complexity of theorem proving procedures". Proceedings of the Third Annual ACM Symposium on Theory of Computing. pp. 151–158. doi:10.1145/800157.805047. ISBN 9781450374644. S2CID 7573663.
\newline
\newline
\textbf{Acknowlegements}: The author wishes to thank RIMS for hosting this conference, Tsunekazu Nishinaka for running it, and Hisa Tsutsui for arranging my attendance.
\end{document}